\documentclass[12pt,reqno]{amsart}

\usepackage[
a4paper,
textwidth=16cm,
textheight=23.3cm,
centering
]{geometry}

\usepackage{amsmath,amssymb,mathtools}
\usepackage{mathrsfs}
\numberwithin{equation}{section}

\usepackage[all,cmtip]{xy}
\usepackage{tikz-cd}
\usepackage{nicematrix}

\usepackage{enumitem}
\usepackage{microtype}
\usepackage{aliascnt}
\usepackage{xcolor}
\usepackage[
colorlinks=true,
linkcolor=red,
citecolor=blue,
urlcolor=blue
]{hyperref}
\usepackage[capitalise,nameinlink]{cleveref}

\theoremstyle{plain}
\newtheorem{theorem}{Theorem}[section]

\newaliascnt{corollary}{theorem}
\newtheorem{corollary}[corollary]{Corollary}
\aliascntresetthe{corollary}

\newaliascnt{lemma}{theorem}
\newtheorem{lemma}[lemma]{Lemma}
\aliascntresetthe{lemma}

\newaliascnt{proposition}{theorem}
\newtheorem{proposition}[proposition]{Proposition}
\aliascntresetthe{proposition}

\theoremstyle{definition}

\newaliascnt{definition}{theorem}
\newtheorem{definition}[definition]{Definition}
\aliascntresetthe{definition}

\newaliascnt{example}{theorem}

\aliascntresetthe{example}

\theoremstyle{remark}

\newaliascnt{remark}{theorem}
\newtheorem{remark}[remark]{Remark}
\aliascntresetthe{remark}

\newaliascnt{observation}{theorem}

\aliascntresetthe{observation}

\crefname{theorem}{Theorem}{Theorems}
\Crefname{theorem}{Theorem}{Theorems}

\crefname{corollary}{Corollary}{Corollaries}
\Crefname{corollary}{Corollary}{Corollaries}

\crefname{lemma}{Lemma}{Lemmas}
\Crefname{lemma}{Lemma}{Lemmas}

\crefname{proposition}{Proposition}{Propositions}
\Crefname{proposition}{Proposition}{Propositions}

\crefname{definition}{Definition}{Definitions}
\Crefname{definition}{Definition}{Definitions}

\crefname{example}{Example}{Examples}
\Crefname{example}{Example}{Examples}

\crefname{remark}{Remark}{Remarks}
\Crefname{remark}{Remark}{Remarks}

\crefname{observation}{Observation}{Observations}
\Crefname{observation}{Observation}{Observations}

\DeclareMathOperator{\Hess}{Hess}
\DeclareMathOperator{\Pic}{Pic}
\DeclareMathOperator{\Cl}{Cl}

\newcommand{\cL}{\mathcal{L}}
\newcommand{\cO}{\mathcal{O}}
\newcommand{\Hmin}{H_{\min}}

\title[Line bundles on regular semisimple Hessenberg varieties]{Line bundles on regular semisimple Hessenberg varieties}
\author[A. Ghosh]{Arkadev Ghosh}
\address{Indian Institute of Technology Delhi,
	Hauz Khas, New Delhi 110016, India}
\email{arkadevghosh198@gmail.com}

\author[P. Saha]{Pinakinath Saha}
\address{Indian Institute of Technology Delhi,
	Hauz Khas, New Delhi 110016, India}
\email{pinakinath@iitd.ac.in}

\subjclass[2020]{Primary 14M15; Secondary 14C22, 14F17}

\keywords{Regular Hessenberg varieties, Picard groups, line bundles,
	flag varieties, vanishing theorems}

\begin{document}

\begin{abstract}
Let $G$ be a connected simply-connected simple algebraic group over $\mathbb{C}$ and  $B$ be a Borel subgroup of $G$. Let $Y$ be an irreducible regular semisimple Hessenberg variety in $G/B$.
In this article, we give a necessary and sufficient condition for the natural restriction map ${\rm res}: \Pic(G/B)\longrightarrow \Pic(Y)$ of Picard groups induced by the inclusion $Y\subseteq G/B$ to be an isomorphism.
\end{abstract}
\maketitle

\section{Introduction}
 
Let $G$ be a connected simply-connected simple algebraic group over $\mathbb{C}$. Let $T$ be a maximal torus of $G$ and $B$ be a Borel subgroup of $G$ containing $T$. Let $\mathfrak{t},\mathfrak{b},\text{and }\mathfrak{g}$ be the Lie algebras of $T$, $B$, and $G$, respectively. 

A vector subspace $H\subseteq \mathfrak{g}$ is said to be a \emph{Hessenberg space} if $\mathfrak{b}\subseteq H$ and ${\rm ad}(\mathfrak{b})(H)\subseteq H$, where ${\rm ad}$ denotes the adjoint representation of $\mathfrak{g}$.
For an element $X$ in $\mathfrak{g}$ and a Hessenberg space $H$, the \emph{Hessenberg variety} associated to $X$ and $H$ is defined as
\[
\Hess(X,H):=\{gB\in G/B:\ \mathrm{Ad}(g^{-1})X\in H\},
\]
where ${\rm Ad}$ denotes the adjoint representation of $G$. 

Assume that $H$ contains every negative simple root space. If $X\in\mathfrak g$ is regular, then $\Hess(X,H)$ is an irreducible projective variety of dimension $\dim(H/\mathfrak b)$, called a \emph{regular Hessenberg variety} (see \cite{Precup}). 

When $X\in\mathfrak t$ is regular semisimple, we call $\Hess(X,H)$ a \emph{regular semisimple Hessenberg variety}. It is a smooth projective $T$-variety under the natural left action of $T$ on $G/B$. Moreover, these varieties have a cell decomposition induced by the cell decomposition of $G/B$ (see \cite{DPS}).  Regular semisimple Hessenberg varieties have attracted much attention, notably through Tymoczko's dot action of the Weyl group on their cohomology \cite{Ty} and its relation with the Stanley--Stembridge and Shareshian--Wachs conjectures (see \cite{BC,GP,SW}).

Assume that $\Hess(X,H)$ is a regular Hessenberg variety.
Since $\Hess(X,H) \subseteq G/B$, there is a natural restriction map
\begin{equation}\label{eq 1.1}
  {\rm res}: \Pic(G/B)\longrightarrow \Pic(\Hess(X, H)).
\end{equation}
If $X\in\mathfrak b$ and $H$ contains every negative simple root space then the Schubert curves are contained in $\Hess(X,H)$ (see \cref{lem:curves}). Then we can easily conclude that the restriction map in \eqref{eq 1.1} is injective (see \cref{lem:restriction-injective}). It is therefore natural to ask when the restriction map in \eqref{eq 1.1} is an isomorphism. 
 
 From now on, we assume that $X\in\mathfrak t$ is regular semisimple and that $\Hess(X,H)$ is irreducible unless stated otherwise.
 
 In type $A$, Ayzenberg, Masuda and Sato \cite{AMS} computed $H^2(\Hess(X,H),\mathbb Z)$ using GKM theory. Their result shows that the restriction map
 \[
 H^2(G/B,\mathbb Z)\longrightarrow H^2(\Hess(X,H),\mathbb Z)
 \]
 is an isomorphism when $H$ contains every negative root space $\mathfrak g_{-\beta}$ with $\operatorname{ht}(\beta)\le2$ (see \cite[Theorem 6.1]{AMS}), and in this case the natural restriction map $$\mathrm{res}: \mathrm{Pic}(G/B)\longrightarrow \mathrm{Pic}(\mathrm{Hess}(X,H))$$ is an isomorphism (see \cite[Corollary 5.3]{AFZ}).

 Our main results are as follows.
 
 \begin{theorem}[\cref{isomorphism:picard group}] \label{main th 1}
 Let $\Hess(X,H)$ be an irreducible regular semisimple Hessenberg variety in $G/B$.
  Then the natural restriction map
 	\[
 	\mathrm{res}:\Pic(G/B)\longrightarrow\Pic(\Hess(X,H))
 	\]
 	of Picard groups is an isomorphism if and only if $\Hmin\subseteq H$, where $\Hmin$ is defined in  \cref{def:Hmin}.
 \end{theorem}
 
 The proof has two parts. First, a root-theoretic argument shows that $\Pic(\Hess(X,H))$ has rank equal to the rank of $G$ if and only if $\Hmin\subseteq H$ (see \cref{rank is n}). Second, under this condition, we show that the intersections of $\Hess(X,H)$ with the Schubert divisors $X(w_0s_i)$ are prime divisors whose classes generate $\Cl(\Hess(X,H))$ (see \cref{irreducibility of intersection}, \cref{class group 2}). We then identify their associated line bundles with the restrictions of $\mathcal L(\omega_i)$, proving surjectivity (see \cref{surjectivity 2}).

Once the Picard groups are identified, it is natural to compare spaces of sections. Brosnan, Hong and Lee \cite{BHL} proved that in type $A$ the restriction
\[
H^0(G/B,\cL(\lambda))\longrightarrow H^0(\Hess(X,H),\cL(\lambda))
\]
is surjective for every regular dominant weight $\lambda$, and $H^{i}(\mathrm{Hess}(X,H),\mathcal{L}(\lambda))=0$ for $i>0$ and for every dominant weight $\lambda$.
Their proof relies on cohomology vanishing of line bundles on Lusztig varieties.
We give an example of an ample line bundle on a Hessenberg variety, when $G$ is of type $G_2$ to illustrate the failure of surjectivity (see \cref{thm:G2}).
    
 
 
 Finally, we prove the following partial result for regular Hessenberg variety $\Hess(X,h)$ in type $A$.
\begin{theorem} [\cref{type A: vanishing}] \label{main theorem 2}
Let $G={\rm SL}(n+1,\mathbb{C})$ $(n\ge 3)$. Let $X\in\mathfrak{g}$ be regular and $h$ be a Hessenberg function of the form
\[
h=(h(1),h(2),h(3),n+1,\dots,n+1),\qquad h(i)\ge i+1, \text{~for~} i= 1,2,3.
\]
Then for any ample line bundle $\cL(\lambda)$ on $G/B$, the restriction map
\[
H^0(G/B,\cL(\lambda))\longrightarrow H^0(\Hess(X,h),\cL(\lambda))
\]
is surjective, and $H^{i}(\Hess(X,h),\cL(\lambda))=0$ for all $i\ge1$.

\end{theorem}

The paper is organized as follows. In \cref{sec:prelim}, we recall the notation and preliminary results. In \cref{sec:Pic}, we prove the Picard-group criterion. In \cref{sec:sections}, we study restriction of sections and higher cohomology, including the example in type $G_2$ and the vanishing theorem in type $A$. We conclude with a question in \cref{sec:concluding}.

\section{Notation and Preliminaries}\label{sec:prelim}
This section recalls some notation and preliminaries. For details on matters related to
algebraic geometry, algebraic groups, Lie algebras, we refer to \cite{Humphreys},\cite{Hum1}, \cite{BK}, \cite{Bou} and \cite{Ha}.

Let $G$, $B$, and $T$ be as in the introduction. Let $W = N_G(T)/T$ be the Weyl group of $G$ with respect to $T$. Denote by $R$ the set of roots of $G$ with respect to $(G, T)$.  Let $R^{+}$ be the set of positive roots with respect to $(B,T)$, and $R^{-}:=-R^{+}$ is the set of negative roots.  Let $\Delta=\{\alpha_1,\dots,\alpha_n\}$ be the set of simple roots in $R^{+}$, where $n$ is the rank of $G$. The simple reflection in $W$ corresponding
to $\alpha_i\in \Delta$ is denoted by $s_i$. Then $W$ is generated by $s_1,\dots , s_n$. 
Let $\ell$ denote the length function, $\le$ the Bruhat order, and $w_0$ the longest element of $W$. 

For each root $\beta\in R$, let $U_\beta$ denote the corresponding root subgroup. Let $P\supseteq B$ be a parabolic subgroup, and set
\[
\Delta_P:=\{\alpha\in \Delta:U_{-\alpha}\subseteq P\}.
\]
Then $W_P:=N_P(T)/T=\langle s_\alpha:\alpha\in \Delta_P\rangle$, and $P$ is generated by $B$ and the root subgroups $U_{-\alpha}$ for $\alpha\in \Delta_P$. For $J\subseteq\Delta$, let $W_J$ be the subgroup of $W$ generated by $\{s_j:\alpha_j\in J\}$, and let $W^J$ be the set of minimal-length left coset representatives of $W/W_J$ in $W$. Denote the longest element of $W_J$ by $w_{0,J}$, and set $w_0^J:=w_0w_{0,J}$. Then $w_{0}^J$ is the minimal length coset representative of $w_0W_J$. 

For $w\in W$, set
\[
R^{+}(w)=\{\beta\in R^{+}:\ w(\beta)\in R^{-}\}.
\]
Then $|R^{+}(w)|=\ell(w)$.

Note that every
$\beta\in  R$ can be expressed uniquely as $\sum_{i=1}^{n}c_i\alpha_i$ such that $c_i\in \mathbb{Z}$ and all of them are
either non-negative or non-positive. This allows us to define the \emph{height} of a root (relative
to $\Delta$) by ${\rm ht}(\beta) = \sum_{i=1}^{n}c_i$.
For $\beta,\gamma\in R^{+}$ we write $\gamma\le\beta$ if $\beta-\gamma$ is a non-negative integer combination of simple roots. Let $\alpha_0$ denote the highest root.

Let $X(T)$ denote the group of all characters of $T$ and set $X(T)_{\mathbb{R}}:=X(T)\otimes_{\mathbb{Z}}\mathbb{R}$.
We identify $X(T)_{\mathbb{R}}$ with $\mathrm{Hom}_{\mathbb{R}}(\mathfrak{t}_{\mathbb{R}},\mathbb{R})$, the dual of the real form of $\mathfrak{t}$. Let $(\cdot,\cdot)$ denote the positive definite $W$-invariant form on $X(T)_{\mathbb R}$ induced by the Killing form on $\mathfrak g$. 

For any $\lambda\in X(T)_{\mathbb{R}}$ and $\beta\in R$, write
\[
\langle\lambda,\beta\rangle:=\frac{2(\lambda,\beta)}{(\beta,\beta)}.
\]
The reflection $s_\beta\in W$ corresponding to $\beta$ acts by
\[
s_\beta(\lambda)=\lambda-\langle\lambda,\beta\rangle\beta.
\]

Let $\omega_1,\dots,\omega_n$ be the fundamental weights, characterized by
\[
\langle\omega_i,\alpha_j\rangle=\delta_{ij}
\qquad (1\le i,j\le n);
\]
where $\delta_{ij}$ denotes the
Kronecker delta.
Let $\rho$ denote the half sum of the positive roots. An integral weight $\lambda\in X(T)$ is called \emph{dominant} (respectively, \emph{regular dominant}) if $\langle\lambda,\alpha_i\rangle\ge0$ (respectively, $\langle\lambda,\alpha_i\rangle>0$) for every $1\le i\le n$. An integral weight $\lambda\in X(T)$ is called \emph{singular} if $\langle\lambda+\rho,\beta\rangle=0$ for some $\beta\in R$, and \emph{non-singular} otherwise.

 For each $w\in W$, choose a representative $\dot w\in N_G(T)$, let $C(w):=B\dot wB/B$ be the corresponding \emph{Schubert cell} in $G/B$. Its closure
 \[
 X(w):=\overline{C(w)}=\bigsqcup_{v\le w}C(v) \subseteq G/B
 \]
 is the \emph{Schubert variety} associated to $w$, of dimension $\ell(w)$. Schubert varieties are irreducible projective and normal, but need not be smooth in general.
 
 In particular, the Schubert varieties $X(w_0s_i)$, $1\le i\le n$, are the \emph{Schubert divisors} of $G/B$, while $X(s_i)\cong\mathbb P^1$ are the \emph{Schubert curves}. 
 
 \subsection{Line bundles on $G/B$}
 Let $\lambda:T\to\mathbb C^\times$ be a character of $T$. Composing $\lambda$ with the canonical projection $B\to T$ gives a character of $B$, also denoted by $\lambda$. Let $\mathbb C_{-\lambda}:=\mathbb{C}$ be the one-dimensional representation of $B$ defined by $b\cdot z=\lambda(b)^{-1}z$. The principal $B$-bundle $G\to G/B$ determines the associated line bundle
 \[
 \mathcal L(\lambda):=G\times^B\mathbb C_{-\lambda},
 \]
 which is the quotient of $G\times\mathbb C$ by the right $B$-action
 \[
 (g,z)\cdot b=(gb,\lambda(b)z).
 \]
 For a subvariety $Y\subseteq G/B$, we write $\mathcal L(\lambda)|_Y$ for its restriction to $Y$. By abuse of notation, we also denote the restriction of $\cL(\lambda)$ to a subvariety of $G/B$ by $\cL(\lambda)$.

 The map $\cL: X(B)\to \Pic (G/B)$ given by $\lambda\mapsto\mathcal L(\lambda)$ is an isomorphism (see \cite[Corollary 3.3, p.~82]{KKLV}). Thus $\Pic(G/B)$ is free abelian with basis $\mathcal L(\omega_1),\dots,\mathcal L(\omega_n)$. Moreover,
 \[
 \mathcal L(\omega_i)\cong\mathcal O_{G/B}(X(w_0s_i)),
 \]
 and
 \begin{equation}\label{eq:degree}
 	\deg\bigl(\mathcal L(\lambda)|_{X(s_i)}\bigr)
 	=\langle\lambda,\alpha_i\rangle
 	\qquad \text{ for all } 1\le i\le n
 \end{equation}
 (see \cite{BK}). 
 
 Under the above isomorphism, dominant weights correspond to globally generated, equivalently nef, line bundles, while regular dominant weights correspond to ample, equivalently very ample, line bundles. Moreover, if $\lambda$ is dominant, then
 \[
 H^0(G/B,\mathcal L(\lambda))\cong V(\lambda)^*,
 \]
 where $V(\lambda)$ is the irreducible representation of $G$ of highest weight $\lambda$.
 
 \subsection{Borel--Weil--Bott}
 We use the Borel--Weil--Bott theorem \cite{Ja} in the following form. Let $\lambda\in X(T)$. If $\lambda$ is singular (i.e., $\langle \lambda+\rho, \beta\rangle=0$ for some $\beta \in R^+$), then
 \[
 H^i(G/B,\cL(\lambda))=0
 \qquad\text{for all }i\ge0.
 \]
 Otherwise, there is a unique $w\in W$ such that $w(\lambda+\rho)-\rho$ is dominant. In this case, \[H^i(G/B,\cL(\lambda))=0 \text{ for } i\ne\ell(w);\]
 and
 \[
 H^{\ell(w)}(G/B,\cL(\lambda))
 \cong V\bigl(w(\lambda+\rho)-\rho\bigr)^*.
 \]
 Moreover,
 \[
 \ell(w)=\bigl|\{\beta\in R^+:\langle\lambda+\rho,\beta\rangle<0\}\bigr|.
 \]
 We call this number the \emph{index} of the non-singular weight $\lambda$. In particular, if $\lambda$ is dominant, then $H^i(G/B,\cL(\lambda))=0$ for all $i>0$.

 \subsection{Hessenberg spaces}

 A subset $I\subseteq R^+$ is a \emph{lower order ideal} if $\beta\in I$ and $\gamma\in R^+$ with $\gamma\le\beta$ imply $\gamma\in I$. For such an ideal, set
 \[
 H_I:=\mathfrak b\oplus\bigoplus_{\beta\in I}\mathfrak g_{-\beta}.
 \]
 The lower ideal condition is equivalent to $H_I$ being stable under the adjoint action of $\mathfrak b$. Moreover, every Hessenberg space is of the form $H_I$ for a unique lower order ideal $I$ (see \cite[Lemma 1]{DPS}). Throughout, we assume that $\Delta\subseteq I$, so that $H_I$ contains every negative simple root space.
 
 For $p\in\mathbb Z_{\ge0}$, define
 \[
 I_p:=\{\beta\in R^+:\operatorname{ht}(\beta)\le p\},
 \qquad H_p:=H_{I_p}.
 \]

 \begin{definition}\label{def:Hmin}
 	Let $I_{\min}$ be the lower order ideal of $R^+$ generated by
 	\[
 	S_{\min}:=\Delta\cup\{s_i(\alpha_j):1\le i\ne j\le n\},
 	\]
 	and set $H_{\min}:=H_{I_{\min}}$.
 \end{definition}
 
 \begin{remark}
 	For $i\ne j$, we have
 	\[
 	s_i(\alpha_j)=\alpha_j-\langle\alpha_j,\alpha_i\rangle\alpha_i.
 	\]
 	If $\alpha_i$ and $\alpha_j$ are orthogonal, then $s_i(\alpha_j)=\alpha_j$. Otherwise, $s_i(\alpha_j)$ is a positive root of height $1-\langle\alpha_j,\alpha_i\rangle\in\{2,3,4\}$. A direct check gives
 	\[
 	H_{\min}=
 	\begin{cases}
 		H_2 & \text{if $G$ is of simply laced type }\\
 		H_2\oplus\mathfrak g_{-(\alpha_{n-1}+2\alpha_n)}
 		& \text{if $G$ is of type } B_n,\\
 		H_2\oplus\mathfrak g_{-(2\alpha_{n-1}+\alpha_n)}
 		& \text{if $G$ is of type } C_n,\\
 		H_2\oplus\mathfrak g_{-(\alpha_2+2\alpha_3)}
 		& \text{if $G$ is of type } F_4,\\
 		H_2\oplus\mathfrak g_{-(2\alpha_1+\alpha_2)}
 		\oplus\mathfrak g_{-(3\alpha_1+\alpha_2)}
 		& \text{if $G$ is of type } G_2.
 	\end{cases}
 	\]
 \end{remark}
 
 \begin{remark}
 	Except in type $G_2$, the set $S_{\min}$ is itself a lower order ideal, so $I_{\min}=S_{\min}$. In type $G_2$, however,
 	\[
 	S_{\min}=\{\alpha_1,\alpha_2,\alpha_1+\alpha_2,3\alpha_1+\alpha_2\}
 	\]
 	is not a lower order ideal: the root $2\alpha_1+\alpha_2$ lies below $3\alpha_1+\alpha_2$ but does not belong to $S_{\min}$.
 \end{remark}
 
 Let $X\in\mathfrak t$ be regular semisimple, and let $I\subseteq R^+$ be a lower order ideal. Then $\Hess(X,H_I)$ is a smooth projective variety of dimension $|I|$ \cite{DPS}. For each $w\in W$, the intersection $C(w)\cap\Hess(X,H_I)$ is isomorphic to an affine space $\mathbb A^{d_w}$, where
 \[
 d_w:=|I\cap R^+(w)|.
 \]
 These intersections form an affine paving of $\Hess(X,H_I)$ \cite{Pr}. As a consequence, its odd Betti numbers vanish, and
 \[
 b_{2k}(\Hess(X,H_I))
 =|\{w\in W:d_w=k\}|.
 \]
 
 If $\Delta\subseteq I$, the first Chern class gives an isomorphism
 \[
 \Pic(\Hess(X,H_I))
 \xrightarrow{\sim}H^2(\Hess(X,H_I),\mathbb Z)
 \]
 by \cite[Corollary 5.3]{AFZ}. Together with the affine paving, this shows that $\Pic(\Hess(X,H_I))$ is free abelian of rank $b_2(\Hess(X,H_I))$.
 
\section{The Picard group}\label{sec:Pic}
In this section, we prove that for an irreducible regular semisimple Hessenberg variety $\Hess(X,H)$, the restriction map
\[
\mathrm{res}:\Pic(G/B)\longrightarrow\Pic(\Hess(X,H))
\]
is an isomorphism if and only if $H_{\min}\subseteq H$.

\begin{lemma}\label{rank is n}
	Let $I\subseteq R^+$ be a lower order ideal containing $\Delta$, and set $H=H_I$. Then:
	\begin{enumerate}
		\item $b_2(\Hess(X,H))\ge n$;
		\item $b_2(\Hess(X,H))=n$ if and only if $H_{\min}\subseteq H$.
	\end{enumerate}
\end{lemma}
\begin{proof}
Proof of (1): For each $1\le i\le n$, we have $R^+(s_i)=\{\alpha_i\}$. Since $\Delta\subseteq I$, it follows that $d_{s_i}=1$. Thus $b_2(\Hess(X,H))\ge n$.

Proof of (2):
Suppose first that $b_2(\Hess(X,H))=n$. For $i\neq j$, we have
	\[
	R^+(s_js_i)=\{\alpha_i,s_i(\alpha_j)\}.
	\]
	If $s_i(\alpha_j)\notin I$, then $d_{s_js_i}=1$. Since $s_js_i$ has length two, this gives an additional one-dimensional cell besides those indexed by the simple reflections, contradicting $b_2(\Hess(X,H))=n$. Hence $s_i(\alpha_j)\in I$ for all $i\ne j$, so $S_{\min}\subseteq I$. As $I$ is a lower order ideal, $I_{\min}\subseteq I$, or equivalently, $H_{\min}\subseteq H$.
	
	Conversely, suppose that $I_{\min}\subseteq I$, and let $w\in W$ satisfy $d_w=1$. Then $w\ne id$, so $w(\alpha_i)\in R^-$ for some simple root $\alpha_i$. Since $\Delta\subseteq I$ and $d_w=1$, this is the unique root in $I$ sent to a negative root by $w$.
	
	Set $u=ws_i$. Then $\ell(u)=\ell(w)-1$ and $u(\alpha_i)=-w(\alpha_i)\in R^+$. We claim that $u=id$. Otherwise, choose a reduced expression
	\[
	u=s_{i_1}\cdots s_{i_k}.
	\]
	Since $u(\alpha_i)\in R^+$, we have $i_k\ne i$. Therefore,
	\[
	s_i(\alpha_{i_k})\in S_{\min}\subseteq I,
	\qquad
	w\bigl(s_i(\alpha_{i_k})\bigr)=u(\alpha_{i_k})\in R^-.
	\]
	Since $i\neq i_k$, the coefficient of $\alpha_{i_k}$ in $s_i(\alpha_{i_k})=\alpha_{i_k}-\langle\alpha_{i_k},\alpha_{i}\rangle\alpha_i$ is $1$. 
	Thus, $\alpha_i$ and $s_i(\alpha_{i_k})$ are two distinct elements of $I\cap R^{+}(w)$. Hence, $d_{w}>1$, which is a contradiction. Hence $u=id$, and $w=s_i$. Thus,
	\[b_2(\mathrm{Hess}(X,H))=|\{w\in W:d_w=1\}|=|\{s_i:\alpha_i\in \Delta\}|=n.\]
\end{proof}

Fix a regular semisimple element $X\in\mathfrak t$. For a Hessenberg space $H$ and $w\in W$, set
\[
Y_w:=C(w)\cap\Hess(X,H).
\]
Since $C(w_0)$ is open in $G/B$, the subset $Y_{w_0}$ is open in $\Hess(X,H)$.

Let $Y:=\Hess(X,H)$. We have
\[
Y\setminus Y_{w_0}
=\bigcup_{i=1}^n\bigl(Y\cap X(w_0s_i)\bigr).
\]
Each Schubert divisor $X(w_0s_i)$ is an effective Cartier divisor on $G/B$. Since $Y_{w_0}\neq \emptyset$ and $Y_{w_0}\cap X(w_0s_i)=\emptyset$, $Y$ is not contained in $X(w_0s_i)$. As $Y$ is an irreducible variety, the restriction of this divisor to $Y$ is an effective Cartier divisor with support $Y\cap X(w_0s_i)$. Thus the boundary $Y\setminus Y_{w_0}$ is pure of codimension one. Let $E_1,\dots,E_m$ be its irreducible components, regarded as prime divisors on $Y$.

\begin{lemma}\label{class group 1}
The class group $\Cl(Y)$ is free abelian with basis given by the classes of $E_1,\dots,E_m$.
\end{lemma}

\begin{proof}
	By \cite[II, Proposition 6.5]{Ha}, there is an exact sequence
	\[
	\bigoplus_{i=1}^m\mathbb Z E_i
	\xrightarrow{f}\Cl(Y)
	\longrightarrow\Cl(Y_{w_0})
	\longrightarrow0,
	\]
	where $f(\sum_i a_iE_i)=[\sum_i a_iE_i]$ and the second map is restriction. Since $Y_{w_0}$ is an affine space, $\Cl(Y_{w_0})=0$. Hence $f$ is surjective.
	
	To prove injectivity, suppose that
	\[
	\sum_{i=1}^m a_iE_i=\operatorname{div}_Y(\varphi)
	\]
	for some rational function $\varphi\in\mathbb C(Y)^\times$. Restricting to $Y_{w_0}$ gives
	\[
	\operatorname{div}_{Y_{w_0}}(\varphi|_{Y_{w_0}})=0.
	\]
	Since $Y_{w_0}$ is normal, $\varphi|_{Y_{w_0}}$ is an invertible regular function. The only invertible regular functions on an affine space are nonzero constants. Thus $\varphi$ is constant on the dense open subset $Y_{w_0}$, and hence is constant as a rational function on $Y$. Therefore,
	\[
	\sum_{i=1}^m a_iE_i=0.
	\]
	Since the group of Weil divisors is free abelian on the prime divisors, all $a_i$ vanish. Hence $f$ is injective.
\end{proof}

\begin{proposition}\label{irreducibility of intersection}
	Let $H$ be a Hessenberg space with $H_{\min}\subseteq H$. Then, for every $1\le i\le n$, the intersection $X(w_0s_i)\cap\Hess(X,H)$ is irreducible, and
	\[
	X(w_0s_i)\cap\Hess(X,H)
	=\overline{C(w_0s_i)\cap\Hess(X,H)}.
	\]
\end{proposition}
\begin{proof} 
	Let $Y:=\Hess(X,H)$ and $d:=\dim Y$.  Fix $1\leq i\leq n$. By the decomposition $X(w_0s_i)=\bigcup_{v\leq w_0s_i}C(v)$, we have
 \begin{equation}\label{irreducibility of intersection:eq 1.1}
 	X(w_0s_i)\cap Y
 	=\bigcup_{v\le w_0s_i}Y_v=\bigcup_{v\le w_0s_i}\overline{Y_v},
 \end{equation}
  Each $\overline{Y_v}$ is irreducible. Let $F_1,\dots,F_l$ denote the irreducible components of $X(w_{0}s_i)\cap \mathrm{Hess}(X,H)$. In view of \cref{irreducibility of intersection:eq 1.1}, for every $1\leq 
 j\leq l$, there exists $v_j\leq w_0s_i$ such that $F_j=\overline{Y_{v_j}}$.

Since $Y_{w_0}\neq \emptyset$ and $Y_{w_0}\cap X(w_0s_i)=\emptyset$, $Y$ is not contained in $X(w_0s_i)$. Hence the intersection $Y\cap X(w_0s_i)$ is the support of an effective Cartier divisor on $Y$, so every irreducible component has dimension $d-1$.
 
On the other hand, $R^{+}(w_0s_i)=R^{+}\setminus\{\alpha_i\}$, so $d_{w_0s_i}=d-1$. Thus the  $n$ cells $Y_{w_0s_1},\dots,Y_{w_0s_n}$ have dimension $d-1$.
 Using Poincar\'e duality $b_2(Y)=b_{2(d-1)}(Y)$. 
 By \cref{rank is n},
 \[\{v\in W:\mathrm{dim}(Y_v)=d-1\}=\{w_0s_i:1\leq i\leq n\}.\]
Thus for every $v<w_0s_i$, $\mathrm{dim}(\overline{Y_v})<d-1$. Hence, $X(w_0s_i)\cap Y$ has a unique irreducible component $\overline{Y_{w_0s_i}}$. 
\end{proof}

For $1\le i\le n$, let $D_i:=X(w_0s_i)\cap\Hess(X,H)$. If $H_{\min}\subseteq H$,  then by \cref{irreducibility of intersection} it follows that  each $D_i$ is irreducible. 

\begin{corollary}\label{class group 2}
	Let $H$ be a Hessenberg space with $H_{\min}\subseteq H$. Then $\Cl(\Hess(X,H))$ is free abelian with basis given by the classes of $D_1,\dots,D_n$.
\end{corollary}

\begin{proof}
Since $G/B\setminus C(w_0)=\cup_{i=1}^{n}X(w_0s_i)$, we have $\mathrm{Hess}(X,H)\setminus Y_{w_0}=\cup_{i=1}^{n}D_i$. From \cref{irreducibility of intersection}, each $D_i$ is a prime divisor of $\mathrm{Hess}(X,H)$. For $i\neq j$, we have $w_0s_iB\in D_i$ and $w_0s_iB\notin D_j$, so the $D_i$'s are pairwise distinct. Hence $D_1,\dots,D_n$ are exactly the prime divisors of 
$\mathrm{Hess}(X,H)$ contained in $\mathrm{Hess}(X,H)\setminus Y_{w_0}$. Now the proof follows from \cref{class group 1}.
\end{proof}
\begin{lemma}\label{lem:curves}
Let $X\in\mathfrak{b}$ and let $H=H_I$ with $\Delta\subseteq I$. Then the Schubert curve $X(s_i)$ is contained in $\Hess(X,H)$ for every $1\le i\le n$.
\end{lemma}
\begin{proof}
	Fix a simple root $\alpha_i$, and choose a representative $\dot s_i\in N_G(T)$ of $s_i$. By the Bruhat decomposition,
	\[
	X(s_i)=B/B\sqcup B\dot s_iB/B.
	\]
	Since $X\in\mathfrak b\subseteq H$, the point $B$ belongs to $\Hess(X,H)$.
	
	Now let $gB\in B\dot s_iB/B$. We may write $gB=b\dot s_iB$ for some $b\in B$. Since $\mathfrak b$ is stable under $\operatorname{Ad}(B)$,
	\[
	\operatorname{Ad}(b^{-1})X\in\mathfrak b.
	\]
	The simple reflection $s_i$ sends $\alpha_i$ to $-\alpha_i$ and every other positive root to a positive root. Hence
	\[
	\operatorname{Ad}(\dot s_i^{-1})\mathfrak b
	=
	\mathfrak t\oplus\mathfrak g_{-\alpha_i}
	\oplus\bigoplus_{\beta\in R^+\setminus\{\alpha_i\}}\mathfrak g_\beta
	\subseteq \mathfrak b\oplus\mathfrak g_{-\alpha_i}
	\subseteq H.
	\]
	Therefore
	\[
	\operatorname{Ad}\bigl((b\dot s_i)^{-1}\bigr)X
	=
	\operatorname{Ad}(\dot s_i^{-1})
	\operatorname{Ad}(b^{-1})X
	\in H.
	\]
	Thus every point of $B\dot s_iB/B$ lies in $\Hess(X,H)$, and consequently $X(s_i)\subseteq\Hess(X,H)$.
\end{proof}

\begin{lemma}\label{lem:restriction-injective}
	Let $X\in\mathfrak b$, and let $H$ be a Hessenberg space containing every negative simple root space. Then the restriction map
	\[
	\mathrm{res}:\Pic(G/B)\longrightarrow\Pic(\Hess(X,H))
	\]
	is injective.
\end{lemma}
\begin{proof}
	Write $\mathcal L=\mathcal L(\sum_{i=1}^n a_i\omega_i)\in\Pic(G/B)$ and suppose that its restriction to $\Hess(X,H)$ is trivial. By \cref{lem:curves}, each Schubert curve $X(s_i)$ lies in $\Hess(X,H)$. Hence
	\[
	a_i=\deg\bigl(\mathcal L|_{X(s_i)}\bigr)=0
	\qquad (1\le i\le n).
	\]
	Thus $\mathcal L$ is trivial, proving injectivity.
\end{proof}

\begin{proposition}\label{surjectivity 2}
	Let $H$ be a Hessenberg space with $H_{\min}\subseteq H$. Then the restriction map
	\[
	\mathrm{res}:\Pic(G/B)\longrightarrow\Pic(\Hess(X,H))
	\]
	is surjective.
\end{proposition}
\begin{proof}
	Set $Y:=\Hess(X,H)$. Since $Y$ is smooth, $\Cl(Y)\cong\Pic(Y)$. By \cref{class group 2}, the line bundles $\mathcal O_Y(D_1),\dots,\mathcal O_Y(D_n)$ form a basis of $\Pic(Y)$. Thus it suffices to show that
	\[
	\mathrm{res}(\mathcal L(\omega_i))\cong\mathcal O_Y(D_i)
	\qquad (1\le i\le n).
	\]
	
	We have $\mathcal L(\omega_i)\cong\mathcal O_{G/B}(X(w_0s_i))$. Since $Y$ is not contained in $X(w_0s_i)$, this Schubert divisor restricts to an effective Cartier divisor $E_i$ on $Y$, and
	\[
	\mathrm{res}(\mathcal L(\omega_i))\cong\mathcal O_Y(E_i).
	\]
	The support of $E_i$ is $D_i$, so $E_i=m_iD_i$ for some integer $m_i\ge1$. By \cref{lem:curves}, the Schubert curve $X(s_i)$ lies in $Y$. Hence \eqref{eq:degree} gives
	\[
	m_i\deg\bigl(\mathcal O_Y(D_i)|_{X(s_i)}\bigr)
	=\deg\bigl(\mathcal L(\omega_i)|_{X(s_i)}\bigr)
	=1.
	\]
	Since the degree is an integer, $m_i=1$. Therefore
	$\mathrm{res}(\mathcal L(\omega_i))\cong\mathcal O_Y(D_i)$ for every $i$. This completes the proof.
\end{proof}

\begin{theorem}\label{isomorphism:picard group}
Let $\Hess(X,H)$ be an irreducible regular semisimple Hessenberg variety in $G/B$. Then the restriction map
\[
\mathrm{res}:\Pic(G/B)\longrightarrow\Pic(\Hess(X,H))
\]
is an isomorphism if and only if $H_{\min}\subseteq H$.
\end{theorem}
\begin{proof}
	Suppose first that the restriction map is an isomorphism. Then
	\[
	b_2(\Hess(X,H))
	=\operatorname{rk}\Pic(\Hess(X,H))
	=\operatorname{rk}\Pic(G/B)
	=n.
	\]
	By \cref{rank is n}, we obtain $H_{\min}\subseteq H$.
	
	Conversely, suppose that $H_{\min}\subseteq H$. By \cref{rank is n}, we have $b_2(\Hess(X,H))=n$. Thus $\Pic(\Hess(X,H))$ and $\Pic(G/B)$ are free abelian groups of the same rank $n$. By \cref{surjectivity 2}, the restriction map is surjective, and hence is an isomorphism.
\end{proof}

\begin{corollary}\label{cor:ample}
Let $X\in\mathfrak{t}$ be regular semisimple. Let $H$ be a Hessenberg space such that $\Hmin\subseteq H$, and let $\lambda$ be a weight. Then
\begin{enumerate}
\item $\cL(\lambda)|_{\Hess(X,H)}$ is nef if and only if it is globally generated, if and only if $\lambda$ is dominant.
\item $\cL(\lambda)|_{\Hess(X,H)}$ is ample if and only if it is very ample, if and only if $\lambda$ is regular dominant.
\end{enumerate}

\end{corollary}
\begin{proof}
	Let $Y:=\Hess(X,H)$. By \cref{lem:curves}, each Schubert curve $X(s_i)$ lies in $Y$, and \eqref{eq:degree} gives
	\[
	\deg\bigl(\cL(\lambda)|_{X(s_i)}\bigr)
	=\langle\lambda,\alpha_i\rangle
	\qquad (1\le i\le n).
	\]
	Proof of (1): If $\lambda$ is dominant, then $\cL(\lambda)$ is globally generated on $G/B$ \cite[II.8.5]{Ja}, so its restriction to $Y$ is globally generated and hence nef. Conversely, if $\cL(\lambda)|_Y$ is nef, then
	\[
	\langle\lambda,\alpha_i\rangle
	=\deg\bigl(\cL(\lambda)|_{X(s_i)}\bigr)\ge0
	\]
	for every $i$. Thus $\lambda$ is dominant.
	
	Proof of (2): If $\lambda$ is regular dominant, then $\cL(\lambda)$ is very ample on $G/B$ \cite[II.8.5]{Ja}, so its restriction to $Y$ is very ample and hence ample. Conversely, if $\cL(\lambda)|_Y$ is ample, then it has positive degree on every curve in $Y$. In particular,
	\[
	\langle\lambda,\alpha_i\rangle
	=\deg\bigl(\cL(\lambda)|_{X(s_i)}\bigr)>0
	\]
	for every $i$. Thus $\lambda$ is regular dominant.
\end{proof} 

\subsection{Rank of $\Pic(\Hess(X,H_1))$}\label{sec:H1}
The variety $\Hess(X,H_1)$ is the smooth projective toric variety associated with the fan of Weyl chambers \cite[Lemma 10, Theorem 11]{DPS}. Its effective acting torus is $T/Z(G)$, and the fan lies in the corresponding real cocharacter space. The rays are spanned by the $W$-translates of the fundamental coweights $\omega_1^\vee,\dots,\omega_n^\vee$. The stabilizer of $\omega_i^\vee$ is the parabolic subgroup $W_{\Delta\setminus\{\alpha_i\}}$. Thus the number of rays is
\[
\sum_{i=1}^n[W:W_{\Delta\setminus\{\alpha_i\}}]
\]
(see \cite[Section 1.1]{BJ}). Since the Picard group of a smooth complete toric variety is free abelian of rank equal to the number of rays minus its dimension \cite[Theorem 4.1.3]{Cox}, we obtain
\[
\operatorname{rk}\Pic(\Hess(X,H_1))
=\sum_{i=1}^n\bigl([W:W_{\Delta\setminus\{\alpha_i\}}]-1\bigr).
\]

\begin{remark}
We have $\mathrm{dim}\  \Hess(X,H_1)=n$, and it is easy to see  
\begin{align*}
	\{w\in W:d_w=n-1\}&=\sqcup_{r=1}^{n}\{uw_{0,\Delta\setminus\{\alpha_r\}}:u\in W^{\Delta\setminus\{\alpha_r\}},
	u\neq w_{0}^{\Delta\setminus\{\alpha_r\}}\}.
\end{align*}
\end{remark}

Finally, the restriction map
\[
\mathrm{res}:\Pic(G/B)\longrightarrow\Pic(\Hess(X,H_1))
\]
is an isomorphism if and only if $n=1$. By \cref{isomorphism:picard group}, this is equivalent to $H_{\min}\subseteq H_1$, or to $s_i(\alpha_j)$ being simple for every $i\ne j$. This occurs precisely when the simple roots are pairwise orthogonal. Since $G$ is simple, this is possible only in rank one, in which case
\[
\Hess(X,H_1)=G/B\cong\mathbb P^1.
\]
 
\section{Cohomology of line bundles}\label{sec:sections}

	Let $\Hess(X,H)$ be an irreducible regular semisimple Hessenberg variety in type $A$, and let $\lambda$ be a regular dominant weight. Brosnan, Hong and Lee \cite{BHL} proved that the restriction map
	\[
	H^0(G/B,\cL(\lambda))
	\longrightarrow H^0(\Hess(X,H),\cL(\lambda)|_{\Hess(X,H)})
	\]
	is surjective. In this section, we show that surjectivity can fail in type $G_2$. We also prove surjectivity and higher-cohomology vanishing for a family of Hessenberg functions in type $A$ and every regular element $X$.
	
	Let $X\in\mathfrak g$ be a regular element, and let $I\subseteq R^+$ be a lower order ideal containing $\Delta$. Choose a non-simple maximal root $\gamma\in I$, and set
	\[
	I':=I\setminus\{\gamma\},\qquad
	Y:=\Hess(X,H_I),\qquad Y':=\Hess(X,H_{I'}).
	\]
	Then $I'$ is a lower order ideal containing $\Delta$. By \cite[Proposition 4.1]{AFZ}, the ideal sheaf of $Y'$ in $Y$ is isomorphic to $\cL(-\gamma)|_Y$. Hence there is an exact sequence
	\begin{equation}\label{eq:AFZ}
		0\longrightarrow\cL(-\gamma)|_Y
		\longrightarrow\cO_Y
		\longrightarrow\iota_*\cO_{Y'}
		\longrightarrow0,
	\end{equation}
	where $\iota:Y'\hookrightarrow Y$ is the closed embedding.
	
	\subsection{Removing the highest root}\label{sec:highest}
	
	Assume that the rank of $G$ is at least two, and let $\alpha_0$ be the highest root. Since $\alpha_0$ is not simple, $R^+\setminus\{\alpha_0\}$ is a lower order ideal containing $\Delta$. Set
	\[
	H_{\alpha_0}:=H_{R^+\setminus\{\alpha_0\}}
	=\mathfrak b\oplus
	\bigoplus_{\beta\in R^+\setminus\{\alpha_0\}}\mathfrak g_{-\beta}.
	\]
	Let $X\in\mathfrak g$ be a regular element, and write $Y:=\Hess(X,H_{\alpha_0})$. Applying \eqref{eq:AFZ} with $I=R^+$ and tensoring with $\cL(\lambda)$ gives, for every weight $\lambda$, an exact sequence of sheaves on $G/B$:
	\begin{equation}\label{eq:highest}
		0\longrightarrow\cL(\lambda-\alpha_0)
		\longrightarrow\cL(\lambda)
		\longrightarrow\iota_*\bigl(\cL(\lambda)|_Y\bigr)
		\longrightarrow0,
	\end{equation}
	where $\iota:Y\hookrightarrow G/B$ is the closed embedding.
	
\begin{proposition}\label{prop:highest}

	Assume that $\operatorname{rk}(G)\ge2$. Let $X\in\mathfrak g$ be a regular element, and let $\lambda\in X(T)$ be a dominant weight. Then:
	\begin{enumerate}
		\item If $G$ is not of type $C_n$, with type $B_2$ regarded as type $C_2$, the restriction map
		\[
		H^0(G/B,\cL(\lambda))
		\longrightarrow H^0(\Hess(X,H_{\alpha_0}),\cL(\lambda)|_{\Hess(X,H_{\alpha_0})})
		\]
		is surjective, and
		\[
		H^q(\Hess(X,H_{\alpha_0}),\cL(\lambda)|_{\Hess(X,H_{\alpha_0})})=0
		\qquad\text{for all }q\ge 1.
		\]
		
		\item If $G$ is of type $C_n$, the same conclusions hold provided that
		\[
		\langle\lambda,\alpha_1\rangle\ge1.
		\]
		In particular, they hold for every regular dominant weight $\lambda$. For $\lambda=\omega_2$, however, the restriction map
		\[
		H^0(G/B,\cL(\omega_2))
		\longrightarrow H^0(\Hess(X,H_{\alpha_0}),\cL(\omega_2)|_{\Hess(X,H_{\alpha_0})})
		\]
		has a one-dimensional cokernel, while
		\[
		H^q(\Hess(X,H_{\alpha_0}),\cL(\omega_2)|_{\Hess(X,H_{\alpha_0})})=0
		\qquad\text{for all }q\ge 1.
		\]
	\end{enumerate}
\end{proposition}

\begin{proof}
	Write $Y:=\Hess(X,H_{\alpha_0})$. Since $\lambda$ is dominant, $H^q(G/B,\cL(\lambda))=0$ for all $q\ge 1$. The long exact sequence associated to \eqref{eq:highest} therefore yields
	\[
	0\to H^0(G/B,\cL(\lambda-\alpha_0))
	\to H^0(G/B,\cL(\lambda))
	\to H^0(Y,\cL(\lambda))
	\to H^1(G/B,\cL(\lambda-\alpha_0))
	\to0.
	\] and 
	\[
	H^q(Y,\cL(\lambda))
	\cong H^{q+1}(G/B,\cL(\lambda-\alpha_0)) \quad \text{for all $q\ge1$} .
	\]

	Thus the cokernel of restriction is $H^1(G/B,\cL(\lambda-\alpha_0))$. To prove the surjectivity and vanishing assertions in (1) and the first part of (2), it therefore suffices to show that
	\[
	H^q(G/B,\cL(\lambda-\alpha_0))=0
	\qquad(q\ge 1).
	\]
	
	By the root tables in \cite{Bou}, outside type $C_n$ we have
	\[
	0\le\langle\alpha_0,\alpha_i\rangle\le1
	\qquad(1\le i\le n).
	\]
	In type $C_n$, we have $\alpha_0=2\omega_1$, so
	\[
	\langle\alpha_0,\alpha_1\rangle=2,
	\qquad
	\langle\alpha_0,\alpha_i\rangle=0
	\quad(i\ne1).
	\]
	Under the hypotheses of (1) or the first part of (2), it follows that
	\[
	\langle\lambda+\rho-\alpha_0,\alpha_i\rangle
	=\langle\lambda,\alpha_i\rangle+1
	-\langle\alpha_0,\alpha_i\rangle
	\ge0
	\qquad(1\le i\le n).
	\]
	If one of these pairings is zero, then $\lambda-\alpha_0$ is singular, so all cohomology of $\cL(\lambda-\alpha_0)$ vanishes by the Borel--Weil--Bott theorem. Otherwise, all these integral pairings are at least one, so $\lambda-\alpha_0$ is dominant and has no higher cohomology. This proves (1) and the first part of (2).
	
	Finally, suppose that $G$ is of type $C_n$ and $\lambda=\omega_2$. Since $\alpha_0=2\omega_1$ and $\alpha_1=2\omega_1-\omega_2$, we have
	\[
	\langle\omega_2+\rho-\alpha_0,\alpha_1\rangle=-1,
	\qquad
	s_1(\omega_2+\rho-\alpha_0)=\rho.
	\]
	Hence the index of $\omega_2-\alpha_0$ is $1$. Therefore Borel--Weil--Bott theorem gives
	\[
	H^1(G/B,\cL(\omega_2-\alpha_0))\cong\mathbb C,
	\] and \[
	H^q(G/B,\cL(\omega_2-\alpha_0))=0,\quad \text{ for all } q\neq 1.
	\]
The exact sequence and isomorphisms above show that restriction has a one-dimensional cokernel and that $H^q(Y,\cL(\omega_2))=0$ for $q\ge 1$.
\end{proof}

\begin{remark}\label{rem:question}
	For a regular semisimple element $X$, we do not know whether the restriction map
	\[
	H^0(G/B,\cL(\lambda))
	\longrightarrow H^0(\Hess(X,H),\cL(\lambda))
	\]
	is surjective for every regular dominant weight $\lambda$ whenever $\Hmin\subseteq H$. In type $A$, surjectivity holds for every $H$ by \cite{BHL}. In type $G_2$, the only Hessenberg space satisfying $\Hmin\subseteq H\subsetneq\mathfrak g$ is $H_4$, for which surjectivity follows from \cref{prop:highest}.
\end{remark}

\subsection{Failure of surjectivity in type $G_2$}
Let $G$ be of type $G_2$, with $\alpha_1$ the short simple root
and $\alpha_2$ the long simple root. Let $X\in\mathfrak g$ be a regular element, and set
\[
H:=H_3
=\mathfrak b\oplus\mathfrak g_{-\alpha_1}
\oplus\mathfrak g_{-\alpha_2}
\oplus\mathfrak g_{-(\alpha_1+\alpha_2)}
\oplus\mathfrak g_{-(2\alpha_1+\alpha_2)}.
\]

With the root numbering in \cite{Hum1}, we have
\[
\alpha_1=2\omega_1-\omega_2,\qquad
\alpha_2=-3\omega_1+2\omega_2,\qquad
\rho=\omega_1+\omega_2.
\]
The positive roots are
\[
\alpha_1,\quad \alpha_2,\quad
\alpha_1+\alpha_2,\quad
2\alpha_1+\alpha_2,\quad
3\alpha_1+\alpha_2,\quad
3\alpha_1+2\alpha_2.
\]
Set
\[
\gamma:=3\alpha_1+\alpha_2=3\omega_1-\omega_2,
\qquad
\alpha_0:=3\alpha_1+2\alpha_2=\omega_2.
\]
Then
\[
I_5=R^+,\qquad
I_4=R^+\setminus\{\alpha_0\},\qquad
I_3=I_4\setminus\{\gamma\}.
\]
Here $\alpha_0$ is the highest root in $R^+$, and $\gamma$ is maximal in $I_4$. 

\begin{proposition}\label{thm:G2}
	With the notation above, we have
	\[
	\dim H^0(G/B,\cL(\rho))=64,
	\qquad
	\dim H^0(\Hess(X,H),\cL(\rho))=70
	\] and $H^q(\Hess(X,H),\cL(\rho))=0$ for all $q\ge 1$.
	In particular, the restriction map
	\[
	H^0(G/B,\cL(\rho))
	\longrightarrow H^0(\Hess(X,H),\cL(\rho))
	\]
	is not surjective. 
\end{proposition}
\begin{proof}
	Write
	\[
	X_1:=\Hess(X,H_4),\qquad
	X_2:=\Hess(X,H_3).
	\]
	Thus $G/B\supset X_1\supset X_2$, with $\dim X_1=5$ and $\dim X_2=4$.
	
	Applying \eqref{eq:AFZ} to these successive inclusions and tensoring with $\cL(\lambda)$ gives
	\begin{align}
		0&\to\cL(\lambda-\alpha_0)
		\to\cL(\lambda)
		\to\iota_*\bigl(\cL(\lambda)|_{X_1}\bigr)
		\to0,
		\label{eq:G2a}\\
		0&\to\cL(\lambda-\gamma)|_{X_1}
		\to\cL(\lambda)|_{X_1}
		\to j_*\bigl(\cL(\lambda)|_{X_2}\bigr)
		\to0,
		\label{eq:G2b}
	\end{align}
	where $\iota:X_1\hookrightarrow G/B$ and $j:X_2\hookrightarrow X_1$ are the closed embeddings. These direct images have the same cohomology as the corresponding sheaves on $X_1$ and $X_2$.
	
	For every weight $\lambda$, \eqref{eq:G2a} gives a long exact sequence
	\begin{equation}\label{eq:G2a-cohomology}
		\begin{aligned}
			0&\to H^0(G/B,\cL(\lambda-\alpha_0))
			\to H^0(G/B,\cL(\lambda))
			\to H^0(X_1,\cL(\lambda))\\
			&\to H^1(G/B,\cL(\lambda-\alpha_0))
			\to H^1(G/B,\cL(\lambda))
			\to H^1(X_1,\cL(\lambda))\\
			&\to H^2(G/B,\cL(\lambda-\alpha_0))
			\to\cdots.
		\end{aligned}
	\end{equation}
	Set
	\[
	\mu:=\rho-\gamma=-2\omega_1+2\omega_2,
	\qquad
	\mu-\alpha_0=-2\omega_1+\omega_2.
	\]
	Taking $\lambda=\rho$ in \eqref{eq:G2b} gives
	\begin{equation}\label{eq:G2b-cohomology}
		\begin{aligned}
			0&\to H^0(X_1,\cL(\mu))
			\to H^0(X_1,\cL(\rho))
			\to H^0(X_2,\cL(\rho))\\
			&\to H^1(X_1,\cL(\mu))
			\to H^1(X_1,\cL(\rho))
			\to H^1(X_2,\cL(\rho))\\
			&\to H^2(X_1,\cL(\mu))
			\to\cdots.
		\end{aligned}
	\end{equation}

\textbf{Cohomology on $G/B$.}
Since $\rho$ and $\rho-\alpha_0=\omega_1$ are dominant,
	\begin{equation}\label{eq:G2-dominant-vanishing}
		H^q(G/B,\cL(\rho))
		=H^q(G/B,\cL(\rho-\alpha_0))=0
		\qquad(q\ge1).
	\end{equation}
	Moreover,
	\[
	s_1(\mu+\rho)=\rho+\omega_2,
	\qquad
	s_1(\mu-\alpha_0+\rho)=\rho.
	\]
	Both weights on the right are regular dominant. By the Borel--Weil--Bott theorem,
	\begin{equation}\label{eq:G2-BWB}
		\begin{aligned}
			H^1(G/B,\cL(\mu))&\cong V(\omega_2)^*,\\
			H^1(G/B,\cL(\mu-\alpha_0))&\cong\mathbb C,\\
			H^q(G/B,\cL(\mu))
			=H^q(G/B,\cL(\mu-\alpha_0))&=0
			\quad  \text{ for all }  q\ne1.
		\end{aligned}
	\end{equation}
	 By the Weyl dimension formula
	 \cite[Section 24.3, Corollary, p.~139]{Hum1},
	\begin{equation}\label{eq:G2-dimensions}
		\dim V(\omega_1)=7,\qquad
		\dim V(\omega_2)=14,\qquad
		\dim V(\rho)=64.
	\end{equation}
	
	\textbf{Step 1: Cohomology of $\cL(\rho)$ on $X_1$.}
	Taking $\lambda=\rho$ in \eqref{eq:G2a-cohomology} and using \eqref{eq:G2-dominant-vanishing} gives
	\[
	0\to V(\omega_1)^*
	\to V(\rho)^*
	\to H^0(X_1,\cL(\rho))
	\to0.
	\]
	Together with \eqref{eq:G2-dimensions}, the same long exact sequence yields
	\begin{equation}\label{eq:G2-X1-rho}
		\dim H^0(X_1,\cL(\rho))=57,
		\qquad
		H^q(X_1,\cL(\rho))=0
		\quad(q\ge 1).
	\end{equation}
	
	\textbf{Step 2: Cohomology of $\cL(\mu)$ on $X_1$.}
	Taking $\lambda=\mu$ in \eqref{eq:G2a-cohomology} and using \eqref{eq:G2-BWB} gives the exact sequence
	\begin{equation}\label{eq:G2-X1-mu-sequence}
		0\to H^0(X_1,\cL(\mu))
		\to\mathbb C
		\xrightarrow{\varphi}V(\omega_2)^*
		\to H^1(X_1,\cL(\mu))
		\to0,
	\end{equation}
	and $H^q(X_1,\cL(\mu))=0$ for $q\ge2$. Let
	$r:=\operatorname{rk}\varphi\in\{0,1\}$. Then
	\begin{equation}\label{eq:G2-X1-mu-dimensions}
		\dim H^0(X_1,\cL(\mu))=1-r,
		\qquad
		\dim H^1(X_1,\cL(\mu))=14-r.
	\end{equation}
	
	\textbf{Step 3: Global sections of $\cL(\rho)$ on $X_2$.}
	Using \eqref{eq:G2-X1-rho}, the long exact sequence
	\eqref{eq:G2b-cohomology} gives
	\[
	\begin{aligned}
		0&\to H^0(X_1,\cL(\mu))
		\to H^0(X_1,\cL(\rho))\\
		&\to H^0(X_2,\cL(\rho))
		\to H^1(X_1,\cL(\mu))
		\to0.
	\end{aligned}
	\]
	and \[ H^q(X_2,\cL(\rho)) \cong H^{q+1}(X_1,\cL(\mu)) \quad (q \ge 1). \] Since $H^q(X_1,\cL(\mu))=0$ for $q\ge2$, we conclude that $H^q(X_2,\cL(\rho))=0$ for all $q\ge 1$.
	By \eqref{eq:G2-X1-rho}, \eqref{eq:G2-X1-mu-dimensions},
	\begin{equation}\label{eq:G2-X2-dimension}
		\dim H^0(X_2,\cL(\rho))
		=57-(1-r)+(14-r)=70.
	\end{equation}
	
	Finally, $\rho$ is regular dominant, so $\cL(\rho)$ is ample on $G/B$ and on $X_2$. By \eqref{eq:G2-dimensions},\eqref{eq:G2-X2-dimension},
	\[
	\dim H^0(G/B,\cL(\rho))=64
	<70=\dim H^0(X_2,\cL(\rho)).
	\]
	Hence the restriction map is not surjective.
\end{proof}

\subsection{Surjectivity and vanishing for special Hessenberg functions}
Let $G={\rm SL}(n+1,\mathbb C)$ with $n\ge3$, and let $X\in\mathfrak g$ be a regular element. 
In this subsection, we prove surjectivity of restriction and higher-cohomology vanishing for regular dominant weights when
\begin{align}\label{special eq 1.1}
h=(h(1),h(2),h(3),n+1,\dots,n+1),
\qquad
h(1)\ge2,\quad h(2)\ge3,\quad h(3)\ge4.
\end{align}

For $1\le i\le j\le n$, set
\[
\alpha_{i,j}:=\alpha_i+\cdots+\alpha_j,
\qquad
R_i:=\{\alpha_{i,j}:i\le j\le n\}.
\]
Thus $\operatorname{ht}(\alpha_{i,j})=j-i+1$. 

A \emph{Hessenberg function} is a weakly increasing map
$h:\{1,\dots,n+1\}\to\{1,\dots,n+1\}$ satisfying
$h(i)\ge i$ for every $i$. The lower order ideal associated to a Hessenberg function $h$ is
\[
I(h):=\{\alpha_{i,j}:1\le i\le j\le n,\ j<h(i)\},
\]
and $\Hess(X,h)=\Hess(X,H_{I(h)})$. Note that, if $h$ is as in \eqref{special eq 1.1}, we have
\begin{equation}\label{eq:typeA-missing-roots}
	R^+\setminus I(h)
	=\{\alpha_{i,j}:1\le i\le3,\ h(i)\le j\le n\}.
\end{equation}
In particular, they are non-simple roots in $R_1\cup R_2\cup R_3$.

\begin{lemma}\label{index 1}
	Let $\lambda\in X(T)$ be a regular dominant weight, and let
	$K=\{\beta_1,\dots,\beta_k\}$ be a set of $k\ge1$ distinct non-simple roots in $R_1\cup R_2\cup R_3$. Then
	\[
	\bigl|\{\beta\in R^+:
	\langle\lambda+\rho-\textstyle\sum_{j=1}^k\beta_j,\beta\rangle<0\}\bigr|
	<k.
	\]
\end{lemma}

\begin{proof}
	Set $\sigma:=\sum_{i=1}^{k}\beta_{i}$. Since $\lambda$ is regular dominant, $\lambda-\rho$ is dominant. Hence
	\[
	\langle\lambda+\rho-\sigma,\beta\rangle
	\ge\langle2\rho-\sigma,\beta\rangle
	\qquad \text{ for all } \beta\in R^+.
	\]
	It therefore suffices to prove the assertion for $\lambda=\rho$. For the rest of the proof, set
	\[
	A:=\{\beta\in R^+:\langle2\rho-\sigma,\beta\rangle<0\},
	\qquad
	A_i:=A\cap R_i\qquad \text{ for all } 1\le i\le n.
	\]
	Since $G$ is simply laced,
	\begin{equation}\label{eq:typeA-negative-pairing}
		\beta\in A
		\quad\text{ if and only if }\quad
		\langle\sigma,\beta\rangle>2\operatorname{ht}(\beta).
	\end{equation}
	
	For $1\le i\le3$, write
	\[
	K_i:=K\cap R_i,\qquad
	\sigma_i:=\sum_{\beta\in K_i}\beta,
	\]
	and set
	\[
	p:=|K_1|,\qquad q:=|K_2|,\qquad r:=|K_3|.
	\]
	Then $\sigma=\sigma_1+\sigma_2+\sigma_3$ and $k=p+q+r$.
	
	Writing $\alpha_{a,b}=e_a-e_{b+1}$, we obtain
	\begin{equation}\label{eq:typeA-root-pairing}
		\langle\alpha_{a,b},\alpha_{m,j}\rangle
		=\delta_{a,m}+\delta_{b+1,j+1}
		-\delta_{a,j+1}-\delta_{b+1,m}
	\end{equation}
	where $\delta_{ij}$ denotes the Kronecker delta.
	Thus we have 
	\begin{equation}\label{eq:typeA-root-pairing1}
		\langle\alpha_{a,b},\alpha_{m,j}\rangle=
		\begin{dcases}
			2 & \text{if }(a,b)=(m,j),\\
			1 & \text{if exactly one of }a=m,\ b=j\text{ holds},\\
			-1 & \text{if }a=j+1\text{ or }m=b+1,\\
			0 & \text{otherwise}.
		\end{dcases}
	\end{equation}
	In particular, if $i\ne m$, each $K_i$ contains at most one
	root whose pairing with $\alpha_{m,j}$ equals $1$.

	\textbf{Roots in $A_1$.}
	Since the roots in $K$ are non-simple, \eqref{eq:typeA-root-pairing1} gives
	\[
	\langle\sigma,\alpha_1\rangle=p-q,
	\qquad
	\langle\sigma,\alpha_{1,2}\rangle\le p+1-r.
	\]
	Thus
	\begin{equation}\label{eq:typeA-row1-special}
		\alpha_1\in A\Longrightarrow p\ge q+3,
		\qquad
		\alpha_{1,2}\in A\Longrightarrow p\ge r+4.
	\end{equation}
	For $\beta=\alpha_{1,j}$ with $j\ge3$, we have
	\[
	\langle\sigma_1,\beta\rangle\le p+1,
	\qquad
	\langle\sigma_2,\beta\rangle\le1,
	\qquad
	\langle\sigma_3,\beta\rangle\le1.
	\]
	Hence $\beta\in A$ implies $2j<p+3$. There are no such indices when $p\le3$, and at most $(p-2)/2$ when $p\ge4$. Consequently,
	\begin{equation}\label{eq:typeA-row1-bound}
		|A_1|
		\le\max\left\{0,\frac{p-2}{2}\right\}
		+|\{\alpha_1,\alpha_{1,2}\}\cap A|.
	\end{equation}
	By \eqref{eq:typeA-row1-special}, we have $A_1=\varnothing$ for $p\le2$ and $|A_1|\le1$ for $p=3$.  For $p\ge4$, \eqref{eq:typeA-row1-bound} gives $|A_1|\le p/2+1\le p-1$. Therefore
	\begin{equation}\label{eq:typeA-row1-deficit}
		|A_1|\le p-1
		\qquad \text{ for all } p\ge1.
	\end{equation}
	
	\textbf{Roots in $A_2$.}
	Similarly,
	\[
	\langle\sigma,\alpha_2\rangle\le1+q-r,
	\qquad
	\langle\sigma,\alpha_{2,3}\rangle\le q+2.
	\]
	Thus
	\begin{equation}\label{eq:typeA-row2-special}
		\alpha_2\in A\Longrightarrow q\ge r+2,
		\qquad
		\alpha_{2,3}\in A\Longrightarrow q\ge3.
	\end{equation}
	For $\beta=\alpha_{2,j}$ with $j\ge4$, we have
	\[
	\langle\sigma_1,\beta\rangle\le1,
	\qquad
	\langle\sigma_2,\beta\rangle\le q+1,
	\qquad
	\langle\sigma_3,\beta\rangle\le1.
	\]
	Hence $\beta\in A$ implies $2(j-1)<q+3$. There are no such indices when $q\le3$, and at most $(q-2)/2$ when $q\ge4$. Consequently,
	\begin{equation}\label{eq:typeA-row2-bound}
		|A_2|
		\le\max\left\{0,\frac{q-2}{2}\right\}
		+|\{\alpha_2,\alpha_{2,3}\}\cap A|.
	\end{equation}
	By \eqref{eq:typeA-row2-special}, we have $A_2=\varnothing$ for $q\le1$, $|A_2|\le1$ for $q=2$, and $|A_2|\le2$ for $q=3$. For $q\ge4$, \eqref{eq:typeA-row2-bound} gives $|A_2|\le q/2+1\le q-1$. Therefore
	\begin{equation}\label{eq:typeA-row2-deficit}
		|A_2|\le q-1
		\qquad \text{ for all } q\ge1.
	\end{equation}
	
	\textbf{Roots in $A_3$.}
	We have
	\[
	\langle\sigma_1,\alpha_3\rangle\le1,
	\qquad
	\langle\sigma_2,\alpha_3\rangle\le1,
	\qquad
	\langle\sigma_3,\alpha_3\rangle=r.
	\]
	Thus $\alpha_3\in A$ implies $r\ge1$. If $\alpha_3\in A$ and $r\le2$, at least one of the first two pairings is positive. By \eqref{eq:typeA-root-pairing}, this gives
	\begin{equation}\label{eq:typeA-row3-exception}
		\alpha_3\in A,\ r\le2
		\quad\Longrightarrow\quad
		\alpha_{1,3}\in K\ \text{or}\ \alpha_{2,3}\in K.
	\end{equation}
	
	For $\beta=\alpha_{3,j}$ with $j\ge4$, we have
	\[
	\langle\sigma_1,\beta\rangle\le1,
	\qquad
	\langle\sigma_2,\beta\rangle\le1,
	\qquad
	\langle\sigma_3,\beta\rangle\le r+1.
	\]
	Hence $\beta\in A$ implies $2(j-2)<r+3$. There are at most $r/2$ such indices. Consequently,
	\begin{equation}\label{eq:typeA-row3-bound}
		|A_3|\le
		\begin{cases}
			r/2,&\text{if }\alpha_3\notin A,\\
			1+r/2,&\text{if }\alpha_3\in A.
		\end{cases}
	\end{equation}
	It follows that $A_3=\varnothing$ if $r=0$. Since $|A_3|$ is an integer, these bounds also give
	\[
	|A_3|\le r-1\qquad(r\ge3).
	\]
	For $r\in\{1,2\}$, we have $|A_3|\le r$, with the stronger bound $|A_3|\le r-1$ whenever $\alpha_3\notin A$.

	\textbf{Roots in $A_m$ ($m\ge4$).}
	Let $\beta=\alpha_{m,j}$ with $m\ge4$. By \eqref{eq:typeA-root-pairing1}, each pairing $\langle\sigma_i,\beta\rangle$, $1\le i\le3$, is at most $1$. If $\beta\in A$, then \eqref{eq:typeA-negative-pairing} gives
	\[
	2\operatorname{ht}(\beta)<\langle\sigma,\beta\rangle\le3.
	\]
	Thus $\operatorname{ht}(\beta)=1$, so $\beta=\alpha_m$. In this case, all three pairings must equal one, which requires
	\[
	\alpha_{1,m}\in K_1,\qquad
	\alpha_{2,m}\in K_2,\qquad
	\alpha_{3,m}\in K_3.
	\]
	Since distinct values of $m$ give distinct roots in each $K_i$, and each $A_m$ contains at most the single root $\alpha_m$, we obtain
	\begin{equation}\label{eq:typeA-later-rows}
		\sum_{m=4}^n|A_m|\le\min(p,q,r).
	\end{equation}
	When $n=3$, the sum on the left is empty.

	Set $c:=\min(p,q,r)$. By \eqref{eq:typeA-later-rows}, it suffices to prove
	\begin{equation}\label{eq:typeA-total-bound}
		|A_1|+|A_2|+|A_3|\le k-1-c.
	\end{equation}
	The preceding estimates give $|A_i|\le|K_i|$ for $i=1,2,3$.
	
	{\bf Case $c=0$.}
	If $p\ge1$, then \eqref{eq:typeA-row1-deficit} gives
	\[
	|A_1|+|A_2|+|A_3|\le(p-1)+q+r=k-1.
	\]
	The same conclusion follows from \eqref{eq:typeA-row2-deficit} if $q\ge1$. Otherwise, $p=q=0$ and $r=k\ge1$. Since neither $\alpha_{1,3}$ nor $\alpha_{2,3}$ belongs to $K$, the exceptional case in \eqref{eq:typeA-row3-exception} cannot occur. Thus $|A_3|\le r-1$, proving \eqref{eq:typeA-total-bound}.
	
	{\bf Case $c=1$.}
	Now $p,q,r\ge1$, so
	\[
	|A_1|+|A_2|+|A_3|
	\le(p-1)+(q-1)+r
	=k-2.
	\]
	Thus \eqref{eq:typeA-total-bound} holds.
	
	{\bf Case $c\ge2$.}
	Set
	\[
	t:=|\{\alpha_1,\alpha_{1,2}\}\cap A|.
	\]
	By \eqref{eq:typeA-row1-special}, we have $p\ge c+2t$. Hence \eqref{eq:typeA-row1-bound} gives
	\[
	|A_1|\le t+\frac{p-2}{2}
	\le p-\frac c2-1.
	\]
	If $\alpha_2\notin A$, then \eqref{eq:typeA-row2-bound} gives
	\[
	|A_2|\le1+\frac{q-2}{2}
	=\frac q2\le q-\frac c2.
	\]
	If $\alpha_2\in A$, then $q\ge r+2\ge c+2$ by \eqref{eq:typeA-row2-special}, and
	\[
	|A_2|\le2+\frac{q-2}{2}
	=\frac q2+1\le q-\frac c2.
	\]
	Finally, since $r\ge c$, \eqref{eq:typeA-row3-bound} gives
	\[
	|A_3|\le1+\frac r2
	\le r-\frac c2+1.
	\]
	Adding these inequalities yields
	\[
	|A_1|+|A_2|+|A_3|
	\le k-\frac{3c}{2}
	\le k-1-c,
	\]
	where the last inequality follows from $c\ge2$.
	
	Thus \eqref{eq:typeA-total-bound} holds in all cases. Together with \eqref{eq:typeA-later-rows}, it gives $|A|\le k-1$, as required.
\end{proof}

\begin{lemma}\label{vanishing:1}
	Let $\lambda\in X(T)$ be a regular dominant weight, and let $\beta_1,\dots,\beta_k$ be distinct non-simple roots in $R_1\cup R_2\cup R_3$, where $k\ge1$. Then
	\begin{equation}\label{eq:typeA-ambient-vanishing}
		H^q\left(G/B,\cL\left(\lambda-\sum_{j=1}^k\beta_j\right)\right)=0
		\qquad \text{ for all } q\ge k.
	\end{equation}
\end{lemma}

\begin{proof}
	Set $\mu:=\lambda-\sum_{j=1}^k\beta_j$. If $\mu$ is singular, all cohomology groups vanish by the Borel--Weil--Bott theorem. Otherwise, the only possible nonzero cohomology occurs in degree
	\[
	q_0:=|\{\beta\in R^+:\langle\mu+\rho,\beta\rangle<0\}|.
	\]
	By \cref{index 1}, we have $q_0<k$, proving \eqref{eq:typeA-ambient-vanishing}.
\end{proof}

\begin{theorem}\label{type A: vanishing}
	Let $G={\rm SL}(n+1,\mathbb C)$ with $n\ge3$, and let $X\in\mathfrak g$ be a regular element. Let
	\[
	h=(h(1),h(2),h(3),n+1,\dots,n+1)
	\]
	be a Hessenberg function satisfying $h(1)\ge2$, $h(2)\ge3$, and $h(3)\ge4$. Then, for every regular dominant weight $\lambda$, the restriction map
	\[
	H^0(G/B,\cL(\lambda))
	\longrightarrow H^0(\Hess(X,h),\cL(\lambda))
	\]
	is surjective, and
	\[
	H^q(\Hess(X,h),\cL(\lambda))=0
	\qquad \text{ for all } q\ge 1.
	\]
\end{theorem}

\begin{proof}
	Fix a regular dominant weight $\lambda$, and let $k_0:=|R^+\setminus I(h)|$. If $k_0=0$, then $\Hess(X,h)=G/B$, and the assertions follow from the Borel--Weil--Bott theorem. Assume that $k_0>0$.
	
	List the missing $k_0$ roots in the following order:
	\[
	(\gamma_1,\dots,\gamma_{k_0})
	=\bigl(\alpha_{1,n},\dots,\alpha_{1,h(1)},\,
	\alpha_{2,n},\dots,\alpha_{2,h(2)},\,
	\alpha_{3,n},\dots,\alpha_{3,h(3)}\bigr).
	\]
	
	For $0\le k\le k_0$, set
	\[ J_0:= R^{+},\qquad
	J_k:=R^+\setminus\{\gamma_1,\dots,\gamma_k\},
	\qquad
	X_k:=\Hess(X,H_{J_k}).
	\]
	A positive root strictly above $\alpha_{i,j}$ has the form $\alpha_{a,b}$ with $a\le i$ and $b\ge j$. If $\alpha_{i,j}$ is missing, then $j\ge h(i)$, and the monotonicity of $h$ gives $b\ge j\ge h(i)\ge h(a)$. Thus $\alpha_{a,b}$ is also missing and occurs earlier in the list. Hence $\gamma_{k+1}$ is maximal in $J_k$ for every $0\le k<k_0$.
	
	It follows that each $J_k$ is a lower order ideal containing $\Delta$. We obtain a chain
	\begin{equation}\label{eq:typeA-filtration}
		G/B=X_0\supset X_1\supset\cdots\supset X_{k_0}
		=\Hess(X,h),
	\end{equation}
	in which each inclusion has codimension one. Applying \eqref{eq:AFZ} and tensoring with $\cL(\nu)$ gives, for every weight $\nu$, an exact sequence
	\begin{equation}\label{eq:onestep}
		0\to\cL(\nu-\gamma_{k+1})|_{X_k}
		\to\cL(\nu)|_{X_k}
		\to(\iota_k)_*\bigl(\cL(\nu)|_{X_{k+1}}\bigr)
		\to0,
	\end{equation}
	where $\iota_k:X_{k+1}\hookrightarrow X_k$ is the closed embedding.

	For every $0\le\ell<k_0$, we claim that
	\begin{equation}\label{eq:typeA-twisted-vanishing}
		H^q\left(X_\ell,\cL\left(\lambda-\sum_{\beta\in S}\beta\right)\right)=0
		\qquad \text{ for all } q\ge|S|
	\end{equation}
	for every nonempty subset
	$S\subseteq\{\gamma_{\ell+1},\dots,\gamma_{k_0}\}$.
	
	We prove this by induction on $\ell$. For $\ell=0$, the claim follows from \cref{vanishing:1}, since the roots in $S$ are distinct non-simple roots in $R_1\cup R_2\cup R_3$.
	
	Suppose that the claim holds for some $\ell<k_0-1$. Let
	$S\subseteq\{\gamma_{\ell+2},\dots,\gamma_{k_0}\}$ be nonempty, and set
	\[
	\nu:=\lambda-\sum_{\beta\in S}\beta.
	\]
	The long exact sequence associated to \eqref{eq:onestep}, with $k=\ell$, contains
	\[
	H^q(X_\ell,\cL(\nu))
	\to H^q(X_{\ell+1},\cL(\nu))
	\to H^{q+1}(X_\ell,\cL(\nu-\gamma_{\ell+1})).
	\]
	For $q\ge|S|$, the first group vanishes by the induction hypothesis applied to $S$. The last group also vanishes, since
	\[
	q+1\ge|S|+1=|S\cup\{\gamma_{\ell+1}\}|.
	\]
	Thus the middle group vanishes, proving the claim.
	
	Taking $\ell=k$ and $S=\{\gamma_{k+1}\}$ in \eqref{eq:typeA-twisted-vanishing}, we obtain
	\begin{equation}\label{eq:typeA-onestep-vanishing}
		H^q(X_k,\cL(\lambda-\gamma_{k+1}))=0
		\qquad \text{ for all } q\ge1,\ 0\le k<k_0.
	\end{equation}
	Taking $\nu=\lambda$ in \eqref{eq:onestep} and using \eqref{eq:typeA-onestep-vanishing}, we obtain surjections
	\[
	H^0(X_k,\cL(\lambda))
	\longrightarrow H^0(X_{k+1},\cL(\lambda))
	\]
	and isomorphisms
	\[
	H^q(X_k,\cL(\lambda))
	\cong H^q(X_{k+1},\cL(\lambda))
	\qquad \text{ for all } q\ge1.
	\]
	Composing these maps along \eqref{eq:typeA-filtration} proves surjectivity of restriction from $G/B$ to $\Hess(X,h)$. It also gives
	\[
	H^q(\Hess(X,h),\cL(\lambda))
	\cong H^q(G/B,\cL(\lambda))
	=0
	\qquad \text{ for all } q\ge1,
	\]
	where the last equality follows from the Borel--Weil--Bott theorem.
\end{proof}
\subsection{Concluding remarks}\label{sec:concluding}
The proof of \cref{type A: vanishing} uses the special form of $h$ only in \cref{index 1}. We therefore ask the following question, which we have not been able to answer.

\noindent\textbf{Question.}
Let $G={\rm SL}(n+1,\mathbb C)$, let $\beta_1,\dots,\beta_k$ be distinct non-simple positive roots, where $k\ge1$, and let $\lambda$ be a regular dominant weight. Is the weight
\[
\lambda-\sum_{j=1}^k\beta_j
\]
either singular or of index less than $k$? In other words, if this weight is non-singular, is
\[
\left|\left\{\beta\in R^+:
\left\langle\lambda+\rho-\sum_{j=1}^k\beta_j,\beta\right\rangle<0
\right\}\right|<k?
\]

By the reduction in the proof of \cref{index 1}, it would suffice to prove the stronger inequality
\[
\left|\left\{\beta\in R^+:
\left\langle2\rho-\sum_{j=1}^k\beta_j,\beta\right\rangle<0
\right\}\right|<k.
\]

The analogous statement fails in types $G_2$ and $F_4$, even when $\lambda=\rho$.

In type $G_2$, with $\alpha_1$ short, take
\[
\gamma:=3\alpha_1+\alpha_2.
\]
The weight $\rho-\gamma=-2\omega_1+2\omega_2$ is non-singular of index is $1$, which is equal to the number of roots subtracted.

In type $F_4$, with the root numbering in \cite{Hum1}, take
\[
\begin{aligned}
	\beta_1&:=\alpha_3+\alpha_4,\\
	\beta_2&:=\alpha_2+2\alpha_3+2\alpha_4,\\
	\beta_3&:=\alpha_1+2\alpha_2+4\alpha_3+2\alpha_4.
\end{aligned}
\]
Their sum is
\[
\sigma:=\beta_1+\beta_2+\beta_3
=-\omega_1-2\omega_2+3\omega_3+3\omega_4.
\]
The weight
\[
\rho-\sigma=2\omega_1+3\omega_2-2\omega_3-2\omega_4
\]
is non-singular and $2\rho-\sigma$ pairs negatively exactly with $\alpha_3$, $\alpha_4$, and $\alpha_3+\alpha_4$. Thus its index is $3$, again equal to the number of roots subtracted.
\section*{Acknowledgements}
A.G. thanks the Indian Institute of Technology Delhi for financial support through the Institute Seed Grant for New Faculty (project MIO3124G SG, under the main SEED Grant project MI02794).

{\bf AI Disclosure.}
The authors used ChatGPT (OpenAI) for assistance with
language editing, LaTeX formatting, and discussion and
clarification of mathematical arguments.

 \bibliographystyle{amsplain}
 \bibliography{references}

@article{AFZ,
	author = {Abe, H. and Fujita, N. and Zeng, H.},
	title = {Geometry of regular {Hessenberg} varieties},
	journal = {Transform. Groups},
	volume = {25},
	year = {2020},
	pages = {305--333}
}

@article{AMS,
	author = {Ayzenberg, A. and Masuda, M. and Sato, T.},
	title = {The second cohomology of regular semisimple {Hessenberg} varieties from {GKM} theory},
	journal = {Proc. Steklov Inst. Math.},
	volume = {317},
	year = {2022},
	pages = {1--20}
}

@book{Bou,
	author = {Bourbaki, N.},
	title = {Lie Groups and Lie Algebras, Chapters 4--6},
	series = {Elements of Mathematics},
	publisher = {Springer},
	address = {Berlin},
	year = {2002}
}

@article{BJ,
	author = {Brion, M. and Joshua, R.},
	title = {Equivariant {Chow} ring and {Chern} classes of wonderful symmetric varieties of minimal rank},
	journal = {Transform. Groups},
	volume = {13},
	year = {2008},
	pages = {471--493}
}

@book{BK,
	author = {Brion, M. and Kumar, S.},
	title = {Frobenius Splitting Methods in Geometry and Representation Theory},
	series = {Progress in Mathematics},
	volume = {231},
	publisher = {Birkh{\"a}user},
	address = {Boston},
	year = {2005}
}

@article{BC,
	author = {Brosnan, P. and Chow, T. Y.},
	title = {Unit interval orders and the dot action on the cohomology of regular semisimple {Hessenberg} varieties},
	journal = {Adv. Math.},
	volume = {329},
	year = {2018},
	pages = {955--1001}
}

@misc{BHL,
	author = {Brosnan, P. and Hong, J. and Lee, D.},
	title = {Geometry of regular semisimple {Lusztig} varieties},
	year = {2025},
	eprint = {2504.15868},
	archivePrefix = {arXiv},
	note = {arXiv:2504.15868},
	url = {https://arxiv.org/abs/2504.15868}
}

@book{Cox,
	author = {Cox, D. A. and Little, J. B. and Schenck, H. K.},
	title = {Toric Varieties},
	series = {Graduate Studies in Mathematics},
	volume = {124},
	publisher = {American Mathematical Society},
	address = {Providence, RI},
	year = {2011}
}

@article{DPS,
	author = {{De Mari}, F. and Procesi, C. and Shayman, M. A.},
	title = {{Hessenberg} varieties},
	journal = {Trans. Amer. Math. Soc.},
	volume = {332},
	year = {1992},
	pages = {529--534}
}

@misc{GP,
	author = {Guay-Paquet, M.},
	title = {A second proof of the {Shareshian--Wachs} conjecture, by way of a new {Hopf} algebra},
	year = {2016},
	eprint = {1601.05498},
	archivePrefix = {arXiv},
	note = {arXiv:1601.05498},
	url = {https://arxiv.org/abs/1601.05498}
}

@book{Ha,
	author = {Hartshorne, R.},
	title = {Algebraic Geometry},
	series = {Graduate Texts in Mathematics},
	volume = {52},
	publisher = {Springer},
	address = {New York},
	year = {1977}
}

@book {Humphreys,
	AUTHOR = {Humphreys, James E.},
	TITLE = {Linear algebraic groups},
	NOTE = {Graduate Texts in Mathematics, No. 21},
	PUBLISHER = {Springer-Verlag, New York-Heidelberg},
	YEAR = {1975},
	PAGES = {xiv+247},
	MRCLASS = {20GXX (14LXX)},
	MRNUMBER = {0396773},
	MRREVIEWER = {T. Ono},
}

@book{Hum1,
	author = {Humphreys, J. E.},
	title = {Introduction to {Lie} Algebras and Representation Theory},
	series = {Graduate Texts in Mathematics},
	volume = {9},
	publisher = {Springer},
	address = {New York},
	year = {1972}
}

@book{Ja,
	author = {Jantzen, J. C.},
	title = {Representations of Algebraic Groups},
	edition = {2},
	series = {Mathematical Surveys and Monographs},
	volume = {107},
	publisher = {American Mathematical Society},
	address = {Providence, RI},
	year = {2003}
}

@incollection{KKLV,
	author = {Knop, F. and Kraft, H. and Vust, T.},
	title = {The {Picard} group of a {$G$}-variety},
	booktitle = {Algebraische Transformationsgruppen und Invariantentheorie},
	series = {DMV Seminar},
	volume = {13},
	publisher = {Birkh{\"a}user},
	address = {Basel},
	year = {1989},
	pages = {77--87}
}

@article{Pr,
	author = {Precup, M.},
	title = {Affine pavings of {Hessenberg} varieties for semisimple groups},
	journal = {Selecta Math. (N.S.)},
	volume = {19},
	year = {2013},
	pages = {903--922}
}

@article{SW,
	author = {Shareshian, J. and Wachs, M. L.},
	title = {Chromatic quasisymmetric functions},
	journal = {Adv. Math.},
	volume = {295},
	year = {2016},
	pages = {497--551}
}

@incollection{Ty,
	author = {Tymoczko, J. S.},
	title = {Permutation actions on equivariant cohomology of flag varieties},
	booktitle = {Toric Topology},
	series = {Contemporary Mathematics},
	volume = {460},
	publisher = {American Mathematical Society},
	address = {Providence, RI},
	year = {2008},
	pages = {365--384}
}

@article{Precup,
	author = {Precup, M.},
	title = {The {Betti} numbers of regular {Hessenberg} varieties are palindromic},
	journal = {Transform. Groups},
	volume = {23},
	number = {2},
	year = {2018},
	pages = {491--499}
}
\end{document}